# Carleman linearization approach for chemical kinetics integration toward quantum computation

Takaki Akiba[1,2], Youhi Morii[1], and Kaoru Maruta[1]

*[1]Tohoku University, Institute of Fluid Science, Sendai 9808577, Japan*

*[2]Tohoku University, School of Engineering, Sendai 9808579, Japan*

ABSTRACT

The Harrow, Hassidim, Lloyd (HHL) algorithm is a quantum algorithm expected to accelerate solving large-scale linear ordinary differential equations (ODEs). To apply the HHL to non-linear problems such as chemical reactions, the system must be linearized. In this study, Carleman linearization was utilized to transform nonlinear first-order ODEs of chemical reactions into linear ODEs. Although this linearization theoretically requires the generation of an infinite matrix, the original nonlinear equations can be reconstructed. For the practical use, the linearized system should be truncated with finite size and analysis precision can be determined by the extent of the truncation. Matrix should be sufficiently large so that the precision is satisfied because quantum computers can treat. Our method was applied to a one-variable nonlinear $\dot{y} = -y^2$ system to investigate the effect of truncation orders in Carleman linearization and time step size on the absolute error. Subsequently, two zero-dimensional homogeneous ignition problems for $H_2$–air and $CH_4$–air gas mixtures were solved. The results revealed that the proposed method could accurately reproduce reference data. Furthermore, an increase in the truncation order in Carleman linearization improved accuracy even with a large time-step size. Thus, our approach can provide accurate numerical simulations rapidly for complex combustion systems.

## 1. INTRODUCTION

Because of the rapidly increasing awareness of global environmental problems, combustion technology requirements have increased considerably. Combustion should be controlled under extreme conditions, such as high pressure, high temperature, and lean-fuel conditions. However, analyzing combustion phenomena under such extreme conditions by using an experimental approach is difficult because of the short characteristic time and nonlinear nature of the phenomena. Therefore, numerical approaches have been proposed for detailed analyses. Numerical analyses of combustion systems under extreme conditions require high-precision methods with detailed information related to chemical reactions. However, detailed analyses of reactive flow systems results in high computational costs because of numerous variables and stiffness of the phenomena. The combustion reaction generally involves 10–1000 chemical species. Furthermore, the dimensional parameters increases the scale of the problems because the thickness of the reaction zones is approximately $10^{-5}$–$10^{-4}$ m, whereas the

Corresponding author: Takaki Akiba

takaki.akiba.q3@dc.tohoku.ac.jp

scale of practical combustion systems is $10^{-1}$–$10^{1}$ m. Because of the thin reaction zone, the required mesh size for reactive flow problems is approximately $30^3$ times finer than that for nonreactive flow problems[1]. The characteristic times of fluid dynamics, molar transport, and chemical reactions differ by approximately $10^{0}$–$10^{-2}$, $10^{-2}$–$10^{-5}$, and $10^{-6}$–$10^{-12}$ s for the fluid dynamics time, molecular-transport time, and chemical reaction time scale, respectively[2].

A possible approach for overcoming these problems is to develop an efficient method or algorithm to evaluate chemical reaction problems without precision loss. As the other approach, powerful machine resources can be used to solve high-cost chemical problems. Rapid developments have been achieved in both quantum computing hardware and software. IBM has shared its roadmap of the scale of quantum computers and will launch a quantum computer with a capacity of more than 1000 qubits in 2023[3]. Although utilization of hardware development is limited, quantum machine resources have been used for high-cost problems. The Harrow, Hassidim, Lloyd (HHL) algorithm is widely known as a powerful solver of large linear equations and expected to be utilized for solving large-scale ordinary differential equations ODEs[4]. The HHL algorithm can be used to solve a large-scale linear problem with $K$ variables within $O(\mathrm{poly}(\log(K)))$ compared with $O(K)$ time required for the best classical algorithm. Studies have focused on specialized reactive flow problems, such as the application to Burgers equations, and primarily pure fluid problems[5,6]. Carleman linearization was used because quantum computers can only be used for linear problems, whereas fluid problems involve nonlinear properties.

Carleman linearization is based on the Taylor expansion and was proposed by Carleman[7]. The Jacobian matrix is another widely considered linearization approach. In the Jacobian approach, only the first derivative of the original system is considered. Therefore, the order considered in the linearization process is limited. Carleman linearization theoretically provides a linearized system with infinite orders, which has the same information as the original nonlinear system, and allows the determination of the degree of orders to which the analysis is considered. Thus, the degree of accuracy considered in the Carleman linearization process can be adjusted based on the nonlinearity of the individual problems and other restrictions. The Carleman linearization method has been widely used in system control[8–10]. Studies have reported the advantage of stability and precision for explicit discretization for time under large time step size conditions under which conventional control methods easily diverge[8]. This advantage of Carleman linearization is suitable for quantum computers because of the frequency of communication between classical computers and quantum computers. Carleman linearization may exhibit a special advantage for a system that has an infinite order or whose variables are correlated with each other. The major problem with the Carleman approach is the large size of the linearized system ($O(N^k)$), where $N$ is the number of variables and $k$ is the truncation order considered in the linearization procedure, as introduced in the next section with equations. When the number of chemical species is 1000, $10^{3k}$ number of variables are required for Carleman linearization. This problem was also referred to in a previous study in which the system control methods are compared for nonlinear systems. However, the development and growth of quantum technologies can solve cost-related problems. Based on these consideration, we applied Carleman linearization to chemical reaction systems, with combustion problems as the final goal. As a trial, a nondimensional reactive system with practical chemical reaction mechanisms, which included multiple elementary reactions and chemical species, was selected in this study.

## 2. FORMULATION

Chemical kinetics are typically first-order ODEs. Therefore, the governing equations for ODEs are summarized as follows:

$$\frac{\mathrm{d}\boldsymbol{u}}{\mathrm{dt}} = \boldsymbol{F}(\boldsymbol{u}), \tag{1}$$

where $\boldsymbol{u}$ is the system variable vector and $F(\boldsymbol{u})$ is a nonlinear function.

In this section, the widely used Jacobian matrix linearization is explained. Next, we explain Carleman linearization. Finally, the differences between the Jacobian and Carleman linearization were explained.

### 2.1. *Linearization using the Jacobian matrix*

The governing equations (equation (1)) can be rewritten with linearization using the Jacobian matrix as follows:

$$\begin{aligned}\frac{\mathrm{d}\boldsymbol{u}}{\mathrm{dt}} &= \boldsymbol{F}(\boldsymbol{u}) \\ &\approx \boldsymbol{F}(\boldsymbol{u}^n) + \left(\frac{\partial \boldsymbol{F}}{\partial \boldsymbol{u}}\right)^n (\boldsymbol{u} - \boldsymbol{u}^n) \\ &= \boldsymbol{F}(\boldsymbol{u}^n) + J^n(\boldsymbol{u} - \boldsymbol{u}^n),\end{aligned} \tag{2}$$

where $J = \frac{\partial \boldsymbol{F}}{\partial \boldsymbol{u}}$ is the Jacobian matrix, and $n$ is the time indicator for the solution vector. By discretizing (2) over time, the differential equation can be transformed as follows:

$$\begin{aligned}&\frac{\boldsymbol{u}^{N+1} - \boldsymbol{u}^N}{\Delta t} = \boldsymbol{F}(\boldsymbol{u}^N) + J^N(\boldsymbol{u}^{N+1} - \boldsymbol{u}^N) \\ &(I - \Delta t J)\boldsymbol{u}^{N+1} = (I - \Delta t J^N)\boldsymbol{u}^N + \Delta t \boldsymbol{F}(\boldsymbol{u}^N), \\ &\boldsymbol{u}^{N+1} = \boldsymbol{u}^N + (I - \Delta t J^N)^{-1}\Delta t \boldsymbol{F}(\boldsymbol{u}^N),\end{aligned} \tag{3}$$

where $N$ and $\Delta t$ represent the time step and the time step size, respectively. We compare equation (3) with other methods.

### 2.2. *Linearization using the Carleman matrix*

Next, we briefly explain the formulation of Carleman linearization. Here, $\boldsymbol{F}(\boldsymbol{x})$ of equation (1) can be transformed as follows:

$$\boldsymbol{F}(\boldsymbol{x}) = \sum A_i \boldsymbol{x}^{\otimes i}, \tag{4}$$

where superscripts of "⊗" represent the Kronecker power, which is expressed as follows:

$$\begin{aligned}&\boldsymbol{x}^{\otimes 1} = \boldsymbol{x}, \\ &\boldsymbol{x}^{\otimes n} = \boldsymbol{x} \otimes \boldsymbol{x}^{\otimes n-1} = \boldsymbol{x}^{\otimes n-1} \otimes \boldsymbol{x}.\end{aligned} \tag{5}$$

Equation (1) can be written down with the expression of equation (5) in linear expression by using Carleman linearization as follows:

$$\frac{\mathrm{d}\boldsymbol{X}}{\mathrm{dt}} = A_c \boldsymbol{X}, \boldsymbol{X} = \begin{pmatrix} \boldsymbol{x} \\ \boldsymbol{x}^{\otimes 2} \\ \vdots \\ \boldsymbol{x}^{\otimes n} \\ \vdots \end{pmatrix}, \tag{6}$$

and

$$A_c = \begin{pmatrix} A_1^1 & A_2^1 & A_3^1 & \cdots & A_n^1 & A_{n+1}^1 & A_{n+2}^1 & \cdots \\ & A_1^2 & A_2^2 & \cdots & A_{n-1}^2 & A_n^2 & A_{n+1}^2 & \cdots \\ & & A_1^3 & & \cdots & & A_n^3 & \cdots \\ & & & \ddots & & & & \vdots \end{pmatrix}, \tag{7}$$

where

$$\begin{aligned}&A_j^1 = A_j \\ &A_j^i = A_j^1 \otimes I^{\otimes i-1} + I \otimes A_j^{i-1} \ (i \geq 2),\end{aligned} \tag{8}$$

for $j \geq 1$. From definition $A_c$, the matrix has infinite rows and columns, and such an infinite matrix cannot be considered in the simulations.

Thus, the matrix should be truncated in the order of $n_t$, and we can summarize the matrix as follows:

$$A_c = \begin{pmatrix} A_1^1 & A_2^1 & A_3^1 & \cdots & A_{n_t}^1 \\ & A_1^2 & A_2^2 & \cdots & A_{n_t-1}^2 \\ & & A_1^3 & & A_{n_t-2}^3 \\ & & & \ddots & \vdots \\ & & & & A_1^{n_t} \end{pmatrix}, \quad (9)$$

$$\boldsymbol{X} = \begin{pmatrix} \boldsymbol{x} \\ \boldsymbol{x}^{\otimes 2} \\ \vdots \\ \boldsymbol{x}^{\otimes n_t} \end{pmatrix}.$$

where $n_t$ is termed as the truncation order hereafter. The number of elements in the matrix is $N_x + N_x^2 + N_x^3 + \cdots + N_x^{n_t} = \frac{N_x^{n_t+1} - N_x}{N_x - 1}$. The final target of this study is the chemical reaction problem. Most elementary reactions involve three chemical species at most. The problem is reduced to an upper triangle matrix as follows:

$$A_c = \begin{pmatrix} A_1^1 & A_2^1 & A_3^1 & & O \\ & A_1^2 & A_2^2 & A_3^2 & \\ & & A_1^3 & A_2^3 & A_3^3 \\ & & & \ddots & \ddots & \ddots \\ O & & & & A_1^{n_t} \end{pmatrix}. \quad (10)$$

When linearization is completed, the system is discretized for time using an explicit approach as follows:

$$\frac{\boldsymbol{X}^{N+1} - \boldsymbol{X}^N}{\Delta t} = A_c^N \boldsymbol{X}^N, \qquad \boldsymbol{X}^{N+1} = (I + \Delta t A_c^N)\boldsymbol{X}^N. \quad (11)$$

Discretization can also be performed using the implicit approach. An $A_c$ can be assumed as a constant through the time between $N$ and $N + 1$ with a small time step size as follows:

$$\frac{\boldsymbol{X}^{N+1} - \boldsymbol{X}^N}{\Delta t} = A_c^{N+1}\boldsymbol{X}^{N+1} \approx A_c^N \boldsymbol{X}^{N+1}, \qquad (I - \Delta t A_c^N)\boldsymbol{X}^{N+1} = \boldsymbol{X}^N, \qquad \boldsymbol{X}^{N+1} = (I - \Delta t A_c^N)^{-1}\boldsymbol{X}^N. \quad (12)$$

By comparing the matrices whose inverse matrix is required for Eqs. (3) and (12), both the Jacobian method and Carleman linearization produce similar expressions: $(I - \Delta t J)$ and $(I - \Delta t A_c^N)$. The differences in the expression is the solution vectors ($\boldsymbol{u}$ and $\boldsymbol{X}$) and matrices used ($J$ and $A_c$). The definitions of $A_c$ and $\boldsymbol{X}$ clearly reveal that Carleman linearization involves higher-order elements in terms of the original solution vector $\boldsymbol{x}$ to be obtained. An implicit approach was used in this study because of the heavy stiffness caused by the chemical reactions mentioned in the following sections. As mentioned, the system size can be estimated $O(N^{n_t})$, and becomes large with a slight increase in the number of variables $N_x$ or the order of truncation $n_t$. The sparse matrix solver in the SciPy library in Python was implemented to solve the implicit approach. The solution vector was obtained by using the direct method because we will utilize the HHL method in a quantum computer in the future, and iterative operations should be eliminated, which could potentially increase the communication between classical and quantum computers. SageMath[11] and its library[12] were used for the Carleman linearization procedure using Python interfaces. Python was used to solve the discretized problem.

# 3. RESULTS AND DISCUSSION

## 3.1. *Nonlinear sample system*

A simple description of the problem is presented first. The target problem is expressed as follows:

$$\frac{\mathrm{d}y}{\mathrm{d}t} = -\alpha y^2$$

$$y(t = 0) = 1.0$$

$$\alpha = 1.0$$

This problem was selected because of the simplicity of the formulation with the smallest nonlinearity

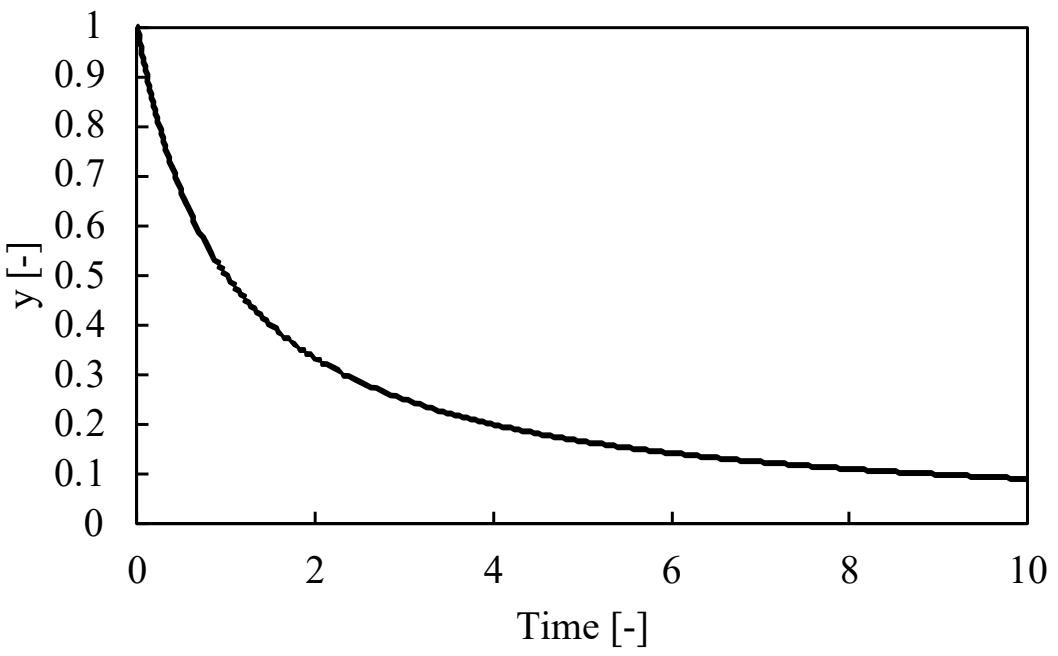

**Figure 2.** Simple nonlinear problem example and solution.

and the consideration of applying the approach to chemically reactive problems. Only single alpha value was applied in this paper because the change of alpha just affects the duration of decay. The initial condition was chosen to simulate the unreacted state of reactants. The solution to this problem is displayed in Fig. 1. In the first step, Carleman linearization was applied to the problem. The Kronecker power of solution $y$ is prepared as follows:

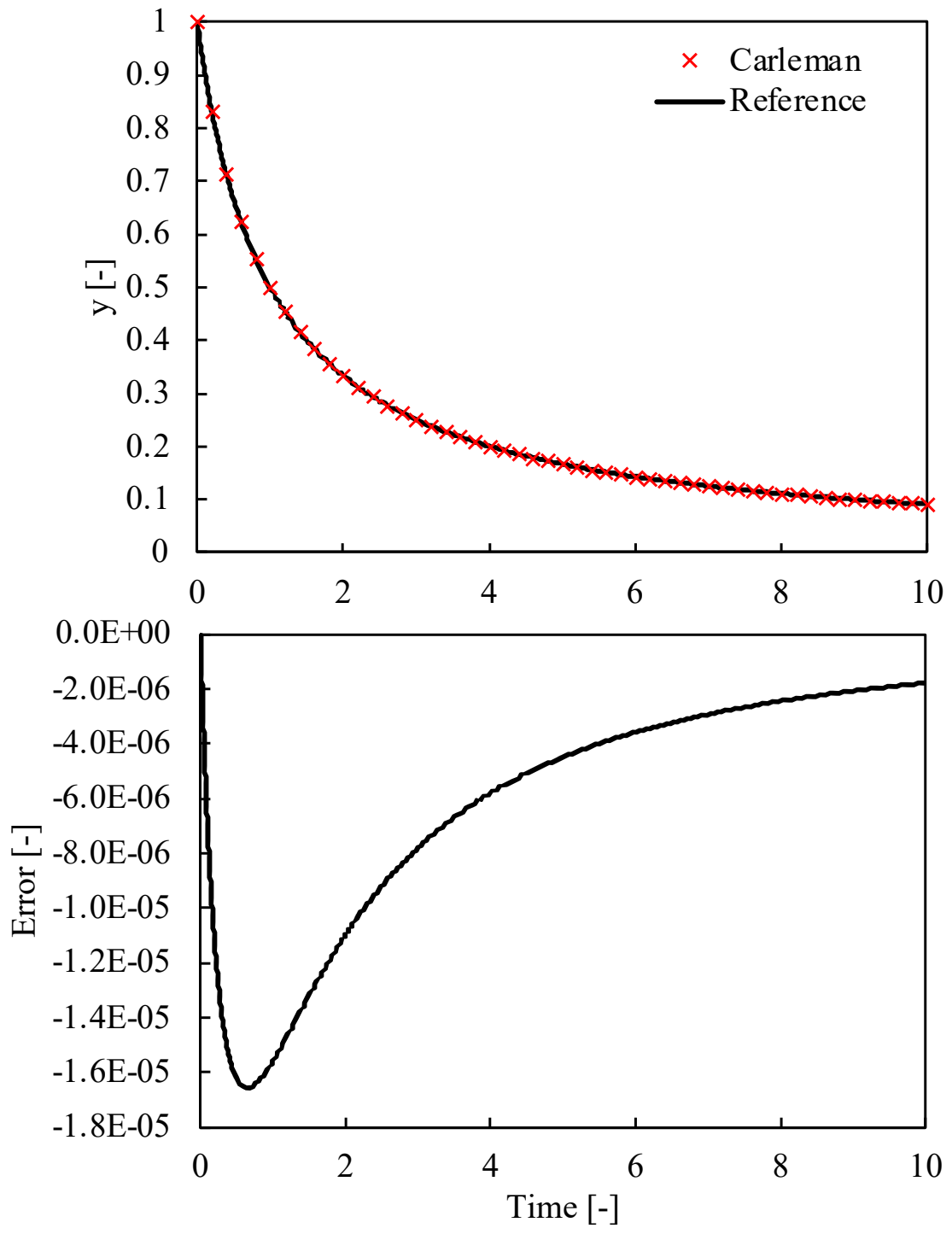

**Figure 3.** Result of Carleman linearization for simple nonlinear problems.

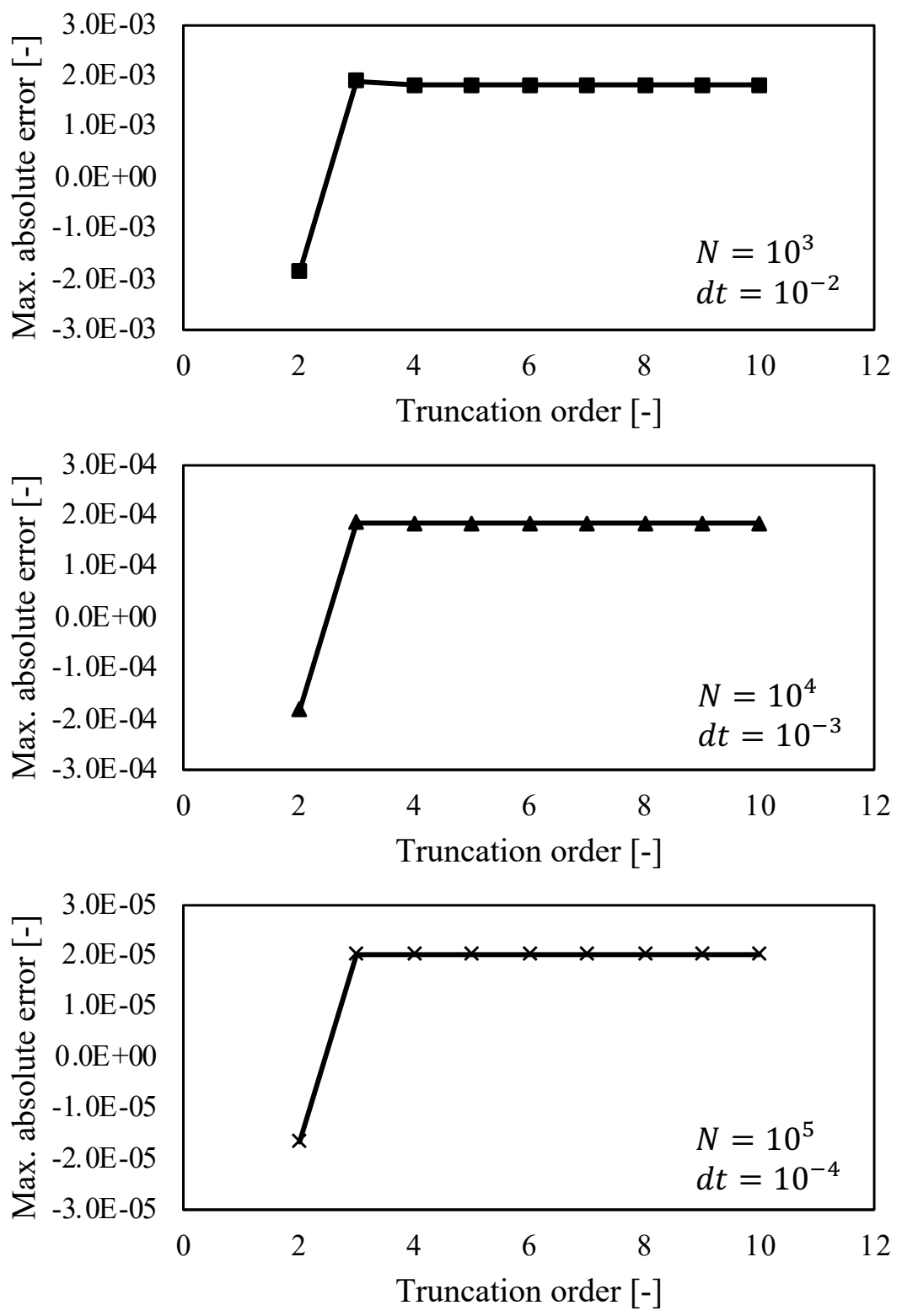

**Figure 1.** Effect of truncation orders and time step sizes on the maximum or minimum absolute value. The representative absolute error is defined as the maximum or minimum absolute error, which has maximum absolute value.

$$\boldsymbol{Y} = (y, y^2, y^3, \cdots, y^n, \cdots)^T.$$

The problem is expressed as follows:

$$\frac{\mathrm{d}\boldsymbol{Y}}{\mathrm{dt}} = A_c \boldsymbol{Y}$$

where

$$A_c = \begin{pmatrix} A_1^1 & A_2^1 & A_3^1 & A_4^1 & \cdots & A_n^1 & A_{n+1}^1 & \cdots \\ & A_1^2 & A_2^2 & A_3^2 & \cdots & A_{n-1}^2 & A_n^2 & \cdots \\ & & \ddots & & & & \vdots & \\ & & & & & A_1^n & A_n^2 & \cdots \\ & & & & & & \ddots & \end{pmatrix}$$

$$= \begin{pmatrix} 0 & -\alpha & 0 & 0 & \cdots & 0 & 0 & \cdots \\ & 0 & -2\alpha & 0 & \cdots & 0 & 0 & \cdots \\ & & \ddots & & & & \vdots & \\ & & & & & 0 & -n\alpha & \cdots \\ & & & & & & \ddots & \end{pmatrix}$$

Because of limited computational resources, the dimensions of the system are truncated in the order of $n$ during computation. The results are compared with reference data in Fig. 2. The reference data were prepared by Euler methods with a time step size of $1.0 \times 10^{-10}$. Fig. 2 displays the absolute errors between the results obtained by Carleman linearization and reference data. The figure clearly displays that the error reaches its maximum at the initial stage of calculation; when $y$ changes considerably. To understand the basic characteristics of the system and the linearization approach, the absolute maximum value of the error was investigated by changing the time step size and truncation order of Carleman linearization.

The results are displayed in Fig. 3. The horizontal axes reveal the truncation order of Carleman linearization. The time step sizes were varied among the three figures. The error decreased with a decrease in the time step size. The error becomes negative for truncation order two, which is the same order as the problem, whereas the error is positive with a truncation order larger than two. With a truncation order of more than two, the absolute maximum errors remained almost constant. Furthermore, the absolute maximum errors were almost the same regardless of the truncation order. These results implied that the truncation order is sufficient with the same order as that of the problem.

### 3.2. *First application to chemical reactive systems - $H_2$–air combustion*

This study is the first to attempt the $H_2$–air combustion problem. The reaction mechanism was based on the USC-II syngas mechanism[13]. To minimize the number of chemical species, only hydrogen, oxygen, and related species were extracted from the base mechanisms. Nitrogen was also added as an inert species in air. The extracted species are displayed in the figure. In actual

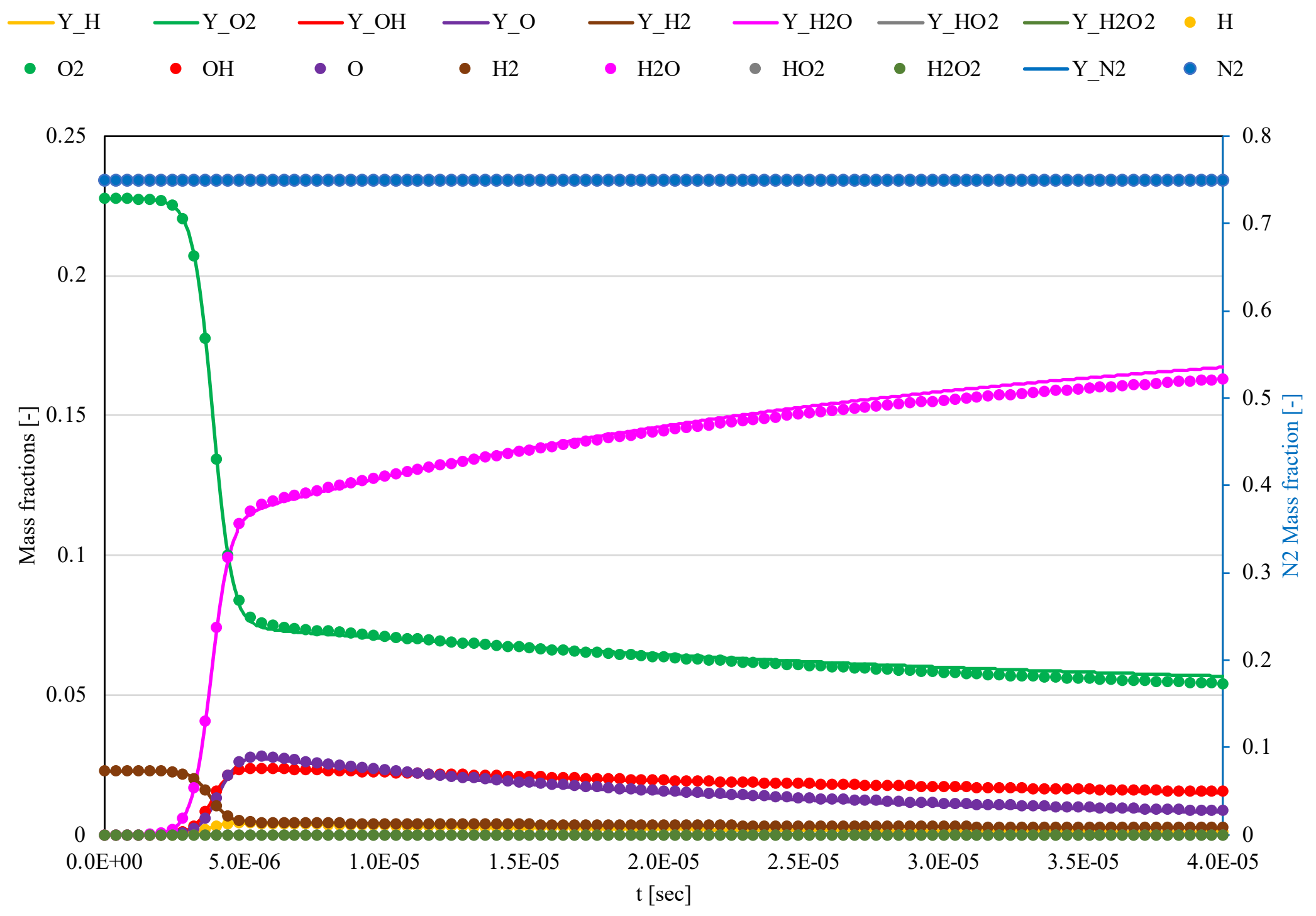

**Figure 4.** Comparison of the results by Carleman linearization and reference for nondimensional iso-thermic transient reactor system of $H_2$–air combustion.

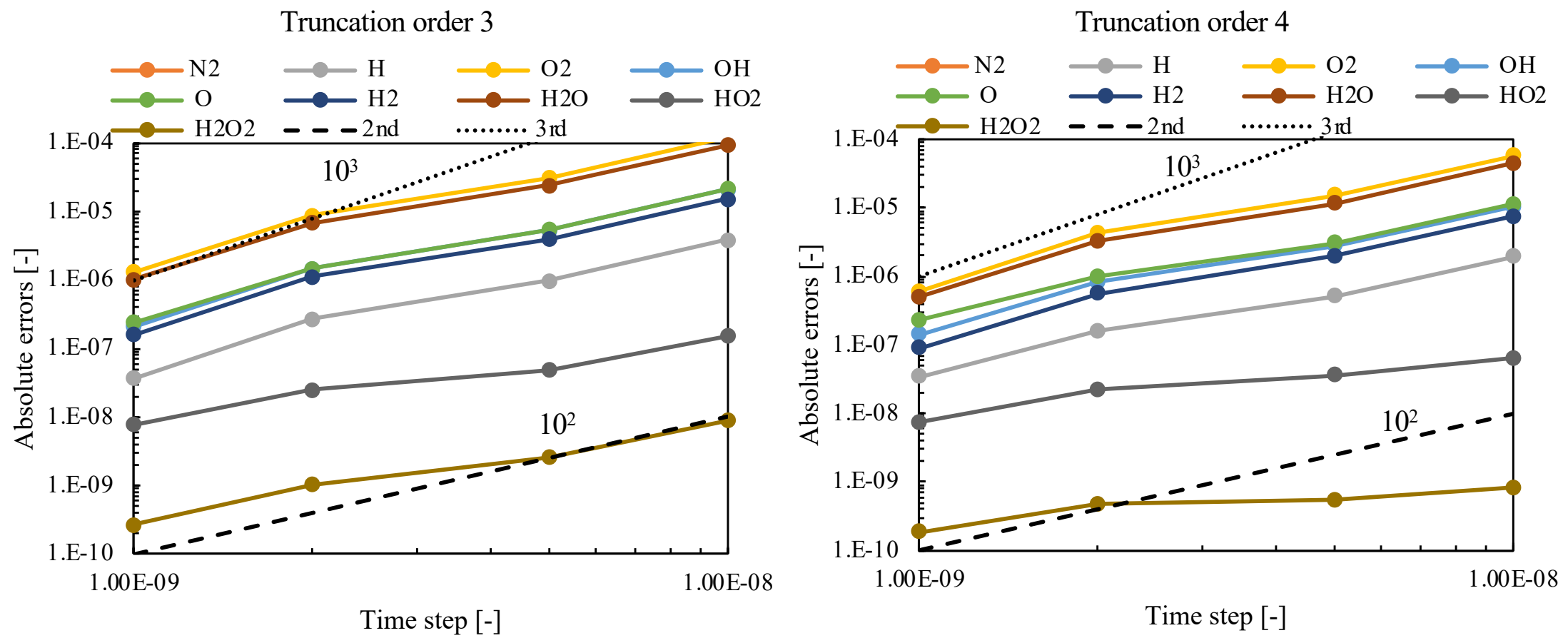

**Figure 5.** Effect of truncation order and time step sizes on the convergence.

chemical reactive systems, the reaction rate of the chemical reaction is expressed as follows:

$$X_A + X_B \leftrightarrow X_C + X_D$$

where $X_i (i = A - D)$ is the chemical species involved and can be evaluated as follows:

$$\begin{aligned}
\dot{\omega_A} &= k_f[X_A][X_B] - k_r[X_C][X_D], \\
\dot{\omega_B} &= k_f[X_A][X_B] - k_r[X_C][X_D], \\
\dot{\omega_C} &= -k_f[X_A][X_B] + k_r[X_C][X_D], \\
\dot{\omega_D} &= -k_f[X_A][X_B] + k_r[X_C][X_D].
\end{aligned} \quad (13)$$

where $\dot{\omega_i}$, $[X_i]$, $k_f$ , and $k_r$ are the reaction rate of chemical species $i$, concentration of chemical species $i$, reaction rate of the forward reaction, and the reaction rate of the reverse reaction, respectively. In actual combustion reaction mechanisms, the number of chemical species involved is in the range of one to three, which indicates that the third order is the maximum order in the evaluation of the chemical reaction. Multiple elementary reactions, such as (13), were considered to investigate the total combustion behavior. As mentioned in the formulation section, we used the implicit method because of its stiffness. Stable solutions were obtained with small time steps. By considering the application to quantum computers in the future, solving a large linear system is preferred to extremely small time steps because quantum machines are powerful for large linear problems, whereas frequent communications are critical. Fig. 4 displays the time history of the chemical species in a zero-dimensional reactor with a time step size of $1.0 \times 10^{-8}$. As displayed in the figure, nine species, including $N_2$ as an inert gas, were selected for the $H_2$–$O_2$ reactions. The equivalence ratio, initial temperature, and pressure were 0.8, 2000 K, and 1 atm, respectively. These conditions are consistent throughout the following discussions, unless otherwise mentioned. Fig. 4 displays the successful evaluation of the chemically reactive system because the overall transition of chemical species and the reaction timing agree.

Similar to the analysis in the previous problem, the effect of the size of the time steps on the simulation error was investigated.

Fig. 5 displays the convergence of the time step sizes for two truncation orders. As references for the slopes, lines 102 and 103 are plotted as dashed and dotted lines, respectively. As expected, the decrease in the time step size improved accuracy. With an increase in the truncation order, the convergence of time step sizes improved considerably, particularly for large time steps. The calculation diverged when the time step increased to more than 2.0 × 108.

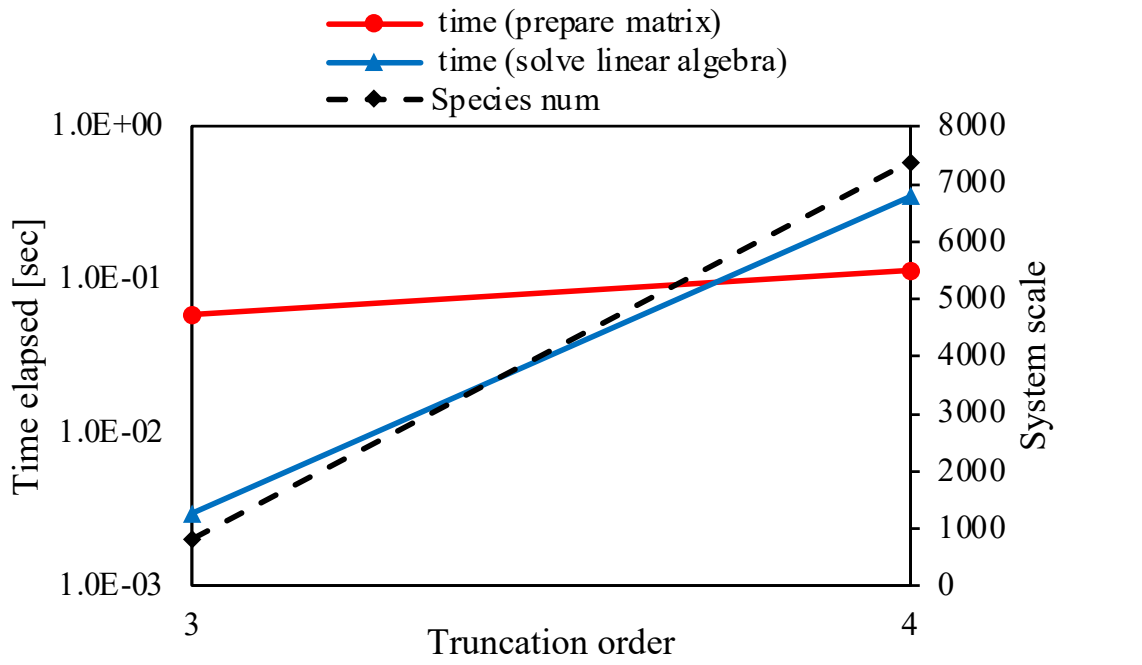

| Truncation | time (prepare matrix) | time (solve linear algebra) | Species num |
|---|---|---|---|
| 3 | 5.91E-02 | 2.97E-03 | 819 |
| | x2 | x120 | x9 |
| 4 | 1.15E-01 | 3.52E-01 | 7380 |

**Figure 6.** Estimation of the effect of truncation order on computational costs.

Although obvious evidence was not obtained, the stiffness of the system influenced the limitation of the timesteps.

Fig. 6 displays the effect of the truncation order on the computational precision and cost. The matrix size is calculated as $(N_v^{n_t+1} - 1)/(N_v - 1)$, where $N_v$ and $n_t$ are the number of variables and the truncation order, respectively, the increase in truncation order by 1 resulted in an increase in the matrix size by a factor of 9. The computational costs increased by a factor of approximately 100, that is, $O(N_v^2)$. An increase in the truncation order had a limited effect on the computational costs for matrix preparation.

### 3.3. *Second application - $CH_4$–air combustion*

The zero-dimensional reactor problem with the $CH_4$–air case was evaluated to confirm the validity of the general hydrocarbon problem and the effect of the number of variables. The chemical reaction model was based on San Diego mech[14]. Only up to a single carbon species was extracted because of the lean-fuel conditions and to minimize computational costs. Thus, 21 species were extracted. The equivalence ratio, initial temperature, and pressure were 0.8, 2000 K, and 1 atm, respectively. The time step size was fixed at $2.0 \times 10^{-8}$, and a further increase in the time step size (e.g. $4.0 \times 10^{-8}$) caused numerical divergence. Fig. 7 displays the time history of the chemical species considered. As displayed in the figure, the result of the Carleman approach was delayed. A possible cause of this delay is the approximation $A_c^{N+1} \approx A_c^N$ in the formulation section around Eq. (12)., which will be investigated in the future. Another possible cause of this delay is the insufficient truncation order of Carleman linearization. However, the order cannot be increased because of machine limitations. This discrepancy will be investigated in the future with the temperature-dependent reaction systems.

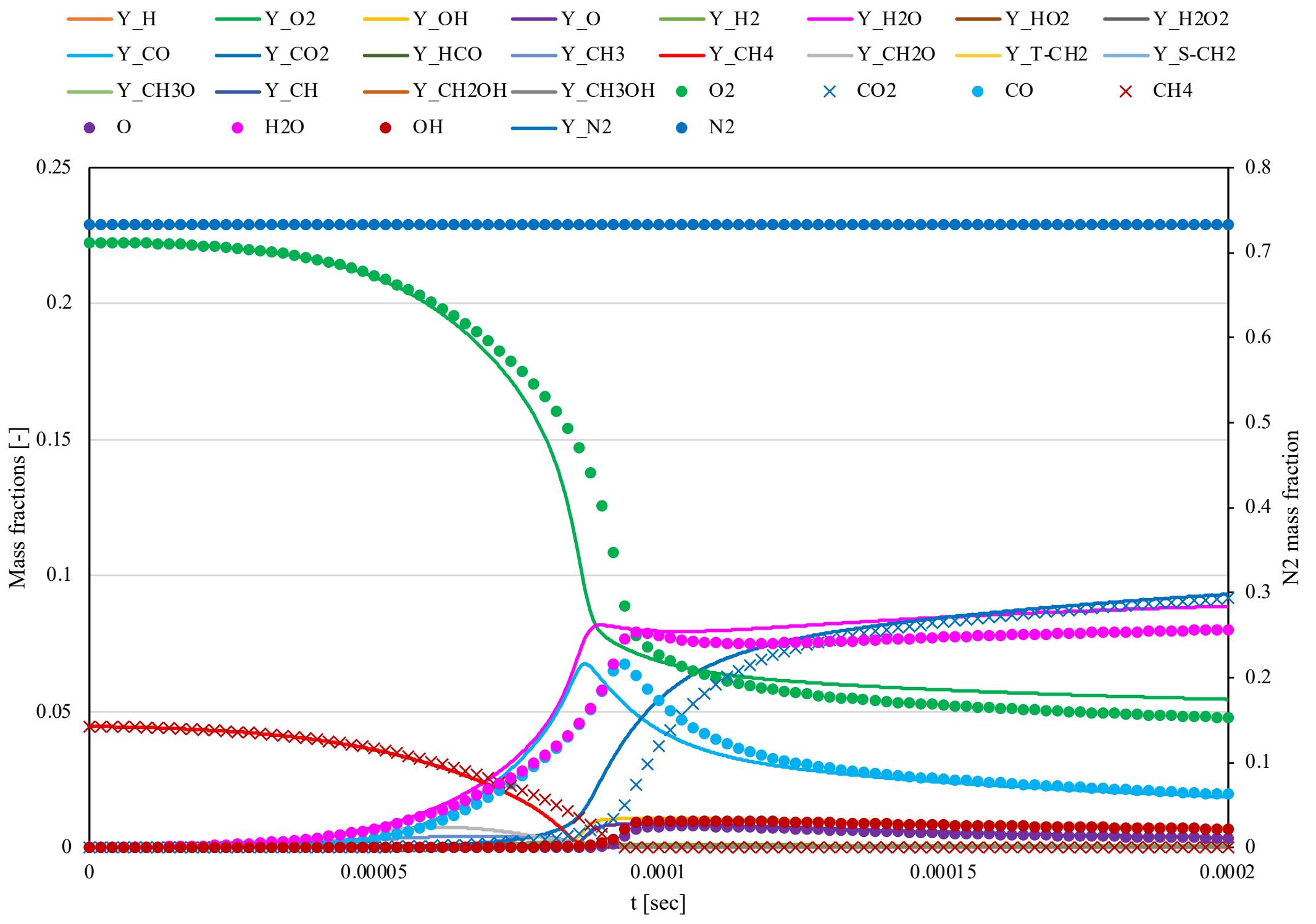

**Figure 7.** Comparison of the results by Carleman linearization and reference for nondimensional isothermic transient reactor system of $CH_4$–air combustion.

## 4. CONCLUSION

In order to utilize the quantum computation resources for large-cost chemical reaction analysis, the nonlinear nature of chemical reactions needs to be linearized for the application of the HHL algorithm, the powerful quantum algorithm for large-scale equations. In this study, the Carleman linearization was applied as a linearization method. The linearization method was evaluated using three simple problems. The results showed the validity and reliability of the proposed approach and implied the potential of using the proposed method in quantum computation of chemical kinetics. In the future, we will focus on evaluating more practical combustion problems.

## Acknowledgments

This study was partially supported by AICE and JSPS KAKENHI (202111887).

## Author contributions

T.A. and Y.M. developed and validated the proposed method. All authors analyzed the results and reviewed the manuscript.

## Data availability

The program data will be available at https://github.com/takakiba/carlin-chem.git

## Competing interests

The authors declare that they have no known competing financial interests or personal relationships that could have appeared to influence the work reported in this paper.